\documentclass[12pt]{article}

\usepackage{a4,graphicx,amsmath,amsfonts,amssymb,bbm}
\usepackage{color,url}

\newcommand{\real}{\mathbbm{R}}
\newcommand{\complex}{\mathbbm{C}}

\newtheorem{theorem}{Theorem}

\newtheorem{lemma}{Lemma}

\begin{document}
\begin{center}
{\bf \Large Model order reduction and sparse orthogonal \\[1ex] expansions 
for random linear dynamical systems}

\vspace{5mm}

{\large Roland~Pulch}

{\small Department of Mathematics and Computer Science, \\
Ernst-Moritz-Arndt-Universit\"at Greifswald, \\
Walther-Rathenau-Str.~47, 17489 Greifswald, Germany. \\
Email: {\tt pulchr@uni-greifswald.de}}

\end{center}

\bigskip\bigskip

%%%%%%%%%%%%%%%%%%%%%%%%%%%%%%%%%%%%%%%%%%%%%%%%%%%%%%%%%%%%%%%%%%%%%%%%%%%%%
%%%                         Abstract                                      %%%
%%%%%%%%%%%%%%%%%%%%%%%%%%%%%%%%%%%%%%%%%%%%%%%%%%%%%%%%%%%%%%%%%%%%%%%%%%%%%

\begin{center}
{Abstract}

\begin{tabular}{p{13cm}}
We consider linear dynamical systems of ordinary differential equations 
or differential algebraic equations. 
Physical parameters are substituted by random variables 
for an uncertainty quantification. 
We expand the state variables as well as a quantity of interest into 
an orthogonal system of basis functions, which depend on the 
random variables. 
For example, polynomial chaos expansions are applicable. 
The stochastic Galerkin method yields a larger linear dynamical system, 
whose solution approximates the unknown coefficients in the expansions.
The Hardy norms of the transfer function provide information about the 
input-output behaviour of the Galerkin system.
We investigate two approaches to construct a sparse representation of the 
quantity of interest, where just a low number of coefficients is non-zero. 
Firstly, a standard basis is reduced by the omission of basis functions, 
whose accompanying Hardy norms are relatively small.
Secondly, a projection-based model order reduction is applied to the 
Galerkin system and allows for the definition of new basis functions 
as a sparse representation.
In both cases, we prove error bounds on the sparse approximation 
with respect to Hardy norms. 
Numerical experiments are demonstrated for a test example 
modelling a linear electric circuit. 

\bigskip

\end{tabular}
\end{center}

\clearpage

%%%%%%%%%%%%%%%%%%%%%%%%%%%%%%%%%%%%%%%%%%%%%%%%%%%%%%%%%%%%%%%%%%%%%%%%%%%%%
%%%                        Introduction                                   %%%
%%%%%%%%%%%%%%%%%%%%%%%%%%%%%%%%%%%%%%%%%%%%%%%%%%%%%%%%%%%%%%%%%%%%%%%%%%%%%

\section{Introduction}
In science and engineering, mathematical modelling often yields 
dynamical systems of ordinary differential equations (ODEs) or 
differential algebraic equations (DAEs). 
We focus on linear time-invariant dynamical systems. 
A quantity of interest is defined as an output of the problem. 
Physical parameters of the systems may exhibit uncertainties 
due to measurement errors or imperfections of an industrial 
manufacturing, for example. 
The uncertainties are described by the introduction of random variables. 
Since many parameters often appear in a system, we are interested in 
the case of high numbers of random variables.

We expand the state variables as well as the quantity of interest 
into an orthonormal system of basis functions depending 
on the random variables. 
For example, the expansions of the polynomial chaos can be used, 
see~\cite{augustinetal,ernst,xiu-karniadakis,xiu-book}. 
Our aim is the construction of a sparse approximation to the random
quantity of interest, where only a few basis functions are required 
for a sufficiently accurate representation. 
Equivalently, we identify a low-dimensional subspace, 
which allows for a good approximation.  
Several previous works exist concerning the computation of such 
a sparse approximation. 
As a tool was used, for example, least angle regression~\cite{blatman}, 
sparse grid quadrature~\cite{conrad-marzouk}, 
compressed sensing~\cite{doostan} 
and $\ell_1$-minimisation~\cite{jakeman2015,jakeman2016}.
Our task can also be seen as an identification of a stochastic 
reduced basis, which was examined for random linear systems of 
algebraic equations in~\cite{nair,sachdeva}.
On the one hand, some methods start from a small set of basis functions 
and extend the basis successively until the approximation becomes 
sufficiently accurate.
On the other hand, some techniques require the choice of an initial set 
of basis functions, which is large and often overprecise, 
and reduce this basis. 
We apply strategies of the latter type.

Either a stochastic Galerkin method or a stochastic collocation technique 
yields approximations to the unknown coefficient functions in the 
expansions, see~\cite{xiu-hesthaven,xiu-book}. 
In this paper, we employ the stochastic Galerkin approach, 
which induces a larger linear dynamical system with many outputs. 
Hardy norms provide a measure for the importance of each output, 
where the $\mathcal{H}_2$-norm and $\mathcal{H}_{\infty}$-norm 
are used. 
Since the system becomes huge for large numbers of random parameters, 
a high potential for a model order reduction (MOR) appears. 
General theory on MOR can be found in the 
monographs~\cite{antoulas,benner,schilders}, for example.
We focus on projection-based techniques for the reduction 
of linear dynamical systems, 
see \cite{feng,freund,gugercin,gugercin-stykel}. 
Projection-based MOR was applied to the stochastic Galerkin system 
in the previous works 
\cite{mi,pulch-maten-augustin,pulch-maten,pulch-scee,zou}.

We investigate two strategies for the construction of a 
sparse approximation. 
Firstly, a large initial basis is reduced by neglecting outputs 
of the Galerkin system with relatively low Hardy norms. 
This reduction implies directly a sparse approximation to the 
random quantity of interest.  
Secondly, a general projection-based MOR technique decreases the 
dimensionality of the Galerkin system.  
We show that this MOR allows for the identification of a sparse approximation 
to the random quantity of interest provided that the reduction achieves 
a sufficiently small system. 
In both cases, error bounds are proved for the sparse approximations 
with respect to Hardy norms.

The paper is organised as follows. 
In Section~\ref{sec:problem-def}, we introduce the problem formulation 
and review already existing theory to be applied. 
The construction of a sparse approximation by omitting basis functions 
is examined in Section~\ref{sec:sparsebasis}. 
The definition of new basis functions using the information from an MOR 
is discussed in Section~\ref{sec:newbasis}.
We present numerical results for an illustrative example 
in Section~\ref{sec:example}.

%%%%%%%%%%%%%%%%%%%%%%%%%%%%%%%%%%%%%%%%%%%%%%%%%%%%%%%%%%%%%%%%%%%%%%%%%%%%%
%%%                     Problem Definition                                %%%
%%%%%%%%%%%%%%%%%%%%%%%%%%%%%%%%%%%%%%%%%%%%%%%%%%%%%%%%%%%%%%%%%%%%%%%%%%%%%

\section{Problem definition}
\label{sec:problem-def}
In this section, we define the problem under investigation. 
Furthermore, results from previous literature, which are relevant 
for our approaches, are outlined.

%%%%%%%%%%%%%%%%%%%%%%%%%%%%%%%%%%%%%%%%%%%%%%%%%%%%%%%%%%%%%%%%%%%%%%%%%%%%%
\subsection{Linear dynamical systems}
We consider a linear time-invariant system in descriptor form
\begin{equation} \label{dae}
\begin{array}{rcl}  
E(p) \dot{x}(t,p) & = & 
A(p) x(t,p) + B(p) u(t) \\[1ex] 
y(t,p) & = & C(p) x(t,p) , \\
\end{array} 
\end{equation}
where the matrices $A,E \in \real^{n \times n}$, 
$B \in \real^{n \times n_{\rm in}}$ and 
$C \in \real^{n_{\rm out} \times n}$ depend on 
physical parameters $p \in \Pi \subseteq \real^q$. 
The input $u : [0,\infty) \rightarrow \real^{n_{\rm in}}$ is supplied, 
while the output is defined by 
$y : [0,\infty) \times \Pi \rightarrow \real^{n_{\rm out}}$.
If the matrix~$E$ is regular, then the system~(\ref{dae}) consists of 
ODEs with state variables 
$x : [0,\infty) \times \Pi \rightarrow \real^n$. 
In our analysis, initial values $x(0,p)=0$ are supposed for all~$p \in \Pi$. 
If the matrix~$E$ is singular, then the system~(\ref{dae}) represents 
DAEs with inner variables~$x$. 
We restrict ourselves to the case of single-input-single-output (SISO) 
with $n_{\rm in} = n_{\rm out} = 1$, 
because generalisations to multiple-input-multiple-output (MIMO) 
are straightforward.
We assume that the matrix pencil $\lambda E(p) - A(p)$ is regular 
for each $p \in \Pi$. 
Moreover, let the system~(\ref{dae}) be stable for all $p \in \Pi$,
i.e., the finite eigenvalues $\Sigma(p) \subset \complex$ 
of the matrix pencil $\lambda E(p) - A(p)$ 
exhibit a negative real part.  

The transfer function 
$H : ( \complex \backslash \Sigma(p)) \rightarrow \complex$ 
characterises the input-output behaviour of the SISO system~(\ref{dae}) 
in the frequency domain, 
see~\cite[Eq.~(4.22)]{antoulas} for explicit ODEs 
or~\cite[Eq.~(2.8)]{freund} for DAEs. 
This transfer function reads as
\begin{equation} \label{transferfcn}
H(s,p) := C(p) \left( s E(p) - A(p) \right)^{-1} B(p) 
\qquad \mbox{for} \;\; s \in \complex \backslash \Sigma(p) . 
\end{equation}
In the case of an SISO system~(\ref{dae}), the transfer function 
becomes a rational function with respect to the frequency~$s$.
%We assume that the system~(\ref{dae}) is proper for each $p \in \Pi$, 
%i.e., the limit of $H$ for $s \rightarrow \infty$ exists.

%%%%%%%%%%%%%%%%%%%%%%%%%%%%%%%%%%%%%%%%%%%%%%%%%%%%%%%%%%%%%%%%%%%%%%%%%%%%%
\subsection{Stochastic modelling and orthogonal expansions}
We substitute the parameters~$p \in \Pi$ of the system~(\ref{dae}) 
by independent random variables $p : \Omega \rightarrow \Pi$ 
on some probability space $(\Omega,\mathcal{A},\mu)$ 
with event space~$\Omega$, sigma-algebra~$\mathcal{A}$ and 
probability measure~$\mu$. 
Let a joint probability density function $\rho : \Pi \rightarrow \real$ 
be given.
For a measurable function $f : \Pi \rightarrow \real$, the expected value 
reads as
$$ \mathbb{E} \left[ f \right] := 
\int_{\Omega} f (p(\omega)) \; \mbox{d}\mu(\omega) = 
\int_{\Pi} f (p) \rho(p) \; \mbox{d}p $$
provided that the integral is finite.
The Hilbert space
$$ \mathcal{L}^2(\Pi,\rho) := \left\{ f : \Pi \rightarrow \real \; : \; 
f \; \mbox{measurable and} \;
\mathbb{E} \left[ f^2 \right] < \infty \right\} $$
is equipped with the inner product
$$ < f , g > \; := \mathbb{E} \left[ f g \right] :=  
\int_{\Pi} f (p) g (p) \rho(p) \; \mbox{d}p 
\qquad \mbox{for} \;\; f,g \in \mathcal{L}^2(\Pi,\rho) . $$
We write the associated norm as
$$ \left\| f \right\|_{\mathcal{L}^2(\Omega)} := \sqrt{ < f , f> } $$
with $\mathcal{L}^2(\Omega)$ as an abbreviation for $\mathcal{L}^2(\Pi,\rho)$.
Concerning the dynamical system~(\ref{dae}), 
we assume that 
$x_1(t,\cdot),\ldots,x_n(t,\cdot), y(t,\cdot) \in \mathcal{L}^2(\Pi,\rho)$
pointwise for $t \in [0,\infty)$. 

Now let a complete orthonormal system 
$( \Phi_i )_{i \in \mathbb{N}} \subset \mathcal{L}^2(\Pi,\rho)$ 
be given. 
It holds that $< \Phi_i , \Phi_j > = 0$ for $i \neq j$ 
and $< \Phi_i , \Phi_j > = 1$ for $i = j$. 
We assume that the first basis function is always the constant 
$\Phi_1 \equiv 1$.
In the polynomial chaos (PC), orthogonal polynomials are applied, 
see~\cite{ernst,xiu-karniadakis,xiu-book}. 
The orthonormal multivariate polynomials are just the products 
of the univariate orthonormal polynomials for each random variable.
It follows that the expansions 
\begin{equation} \label{pce} 
x(t,p) = \sum_{i=1}^\infty v_i(t) \Phi_i(p) 
\qquad \mbox{and} \qquad 
y(t,p) = \sum_{i=1}^\infty w_i(t) \Phi_i(p) ,
\end{equation}
where the coefficient functions 
$v_i : [0,\infty) \rightarrow \real^n$ and 
$w_i : [0,\infty) \rightarrow \real$ are defined by
\begin{equation} \label{coefficients}
v_{i,j}(t) = \; < x_j(t,\cdot) , \Phi_i(\cdot) > 
\qquad \mbox{and} \qquad 
w_i(t) = \; < y(t,\cdot) , \Phi_i(\cdot) > , 
\end{equation}
converge in $\mathcal{L}^2(\Pi,\rho)$ pointwise for $t \in [0,\infty)$ 
and $j=1,\ldots,n$. 

%%%%%%%%%%%%%%%%%%%%%%%%%%%%%%%%%%%%%%%%%%%%%%%%%%%%%%%%%%%%%%%%%%%%%%%%%%%%%
\subsection{Sparse orthogonal representations}
Concerning the series~(\ref{pce}), finite approximations result by a 
truncation. 
We employ the approximations 
\begin{equation} \label{pce-appr} 
x^{(\mathcal{I})}(t,p) = \sum_{i \in \mathcal{I}} v_i(t) \Phi_i(p) 
\qquad \mbox{and} \qquad 
y^{(\mathcal{I})}(t,p) = \sum_{i \in \mathcal{I}} w_i(t) \Phi_i(p) 
\end{equation}
with a finite set of indices $\mathcal{I} \subseteq \mathbbm{N}$.
We focus on the output of the random dynamical system~(\ref{dae}) 
as quantity of interest. 
The associated truncation error reads as
\begin{equation} \label{truncationerror}
\left\| y(t,\cdot) - y^{(\mathcal{I})}(t,\cdot) 
\right\|_{\mathcal{L}^2(\Omega)} = 
\sqrt{\sum_{i \notin \mathcal{I}} w_i(t)^2} 
\end{equation}
for each~$t \ge 0$ due to the orthonormality of the basis functions. 

We consider some set~$\mathcal{I}$, which is assumed to generate
a highly accurate approximation to the output. 
If multivariate polynomials are applied in the expansion, then often all 
polynomials up to a total degree~$d$ are included, 
which yields the index set, cf.~\cite[Sect.~2]{blatman},
\begin{equation} \label{index-large}
\mathcal{I}_d := \left\{ i \; : \; 
\Phi_i(p) = \phi_{j_1}^{(1)}(p_1) \cdots \phi_{j_q}^{(q)}(p_q) 
\;\; \mbox{with} \;\; j_1 + \cdots + j_q \le d 
\right\} . 
\end{equation}
Therein, the systems $(\phi_j^{(\ell)})_{j \in \mathbbm{N}_0}$ include 
the univariate orthonormal polynomials with respect to the $\ell$th 
random variable and the degree of $\phi_j^{(\ell)}$ is exactly~$j$.
The number of basis polynomials becomes,
see~\cite[Eq.~(5.24)]{xiu-book},
$$ | \mathcal{I}_d | = \frac{(q+d)!}{q! d!} . $$
Thus the number of basis functions is large for high numbers~$q$ 
of random variables even if the total degree is low, say $d \le 5$.

We start from a (usually large) initial set~$\mathcal{I}$ 
in~(\ref{pce-appr}), where the truncation error~(\ref{truncationerror}) 
is below some given threshold~$\delta > 0$. 
For polynomial expansions, the set $\mathcal{I}_d$ with the smallest 
total degree~$d$ satisfying this accuracy can be chosen.
Our aim is to identify a sparse approximation by one of the following 
two strategies. 
\begin{enumerate}
\item
For high dimensions~$q$, often an index set 
$\mathcal{I}' \subset \mathcal{I}$ exists with 
$| \mathcal{I}' | \ll | \mathcal{I} |$, while the  
truncation error~(\ref{truncationerror}) is still lower than $\delta$ 
for $y^{(\mathcal{I}')}$. 
Since the new approximation can be written in the form
$$ y^{(\mathcal{I}')} (t,p) = 
\sum_{i \in \mathcal{I}} \tilde{w}_i(t) \Phi_i(p) 
\quad \mbox{with} \quad \tilde{w}_i = \left\{
\begin{array}{cl} 
w_i & \mbox{for}\; i \in \mathcal{I}' , \\
0 & \mbox{for}\; i \notin \mathcal{I}' , \\
\end{array}
\right. $$
where most of the coefficients are equal to zero,  
$y^{(\mathcal{I}')}$ is called a sparse representation.
A measure for the sparsity is the ratio $\sigma \in (0,1)$ defined by 
$\sigma := | \mathcal{I}' | / | \mathcal{I} |$,
%$$ \sigma := \frac{| \mathcal{I}' | }{ | \mathcal{I}_d | } , $$ 
see~\cite[Eq.~(30)]{blatman}.
\item
A new finite-dimensional subspace spanned by the 
multivariate orthonormal polynomials 
$\{ \Psi_1 , \ldots , \Psi_r \}$ is constructed, i.e., 
$$ {\rm span} \{ \Psi_1 , \ldots , \Psi_r \} \subset 
{\rm span} \{ \Phi_i \, : \, i \in \mathcal{I} \} . $$
The alternative approximation reads as 
$$ y^{(r)} (t,p) = 
\sum_{i=1}^r \tilde{w}_i(t) \Psi_i(p) $$
with its own coefficient functions $\tilde{w}_1,\ldots,\tilde{w}_r$ 
to be defined. 
The aim is to keep the dimension~$r$ as small as possible such that 
the $\mathcal{L}^2(\Omega)$-norm of the error for this approximation 
is still lower than $\delta$. 
Again the sparsity is measured by $\sigma := r / | \mathcal{I} |$, 
if the initial set $\mathcal{I}$ is used to construct the approximation.
\end{enumerate}
Since the first strategy can be seen as a special case of the 
second strategy, where the functions $\{ \Psi_1 , \ldots , \Psi_r \}$ 
represent just a subset of $\{ \Phi_i \, : \, i \in \mathcal{I} \}$, 
we expect a higher potential for a dimension reduction of the 
subspaces in the second approach.

%%%%%%%%%%%%%%%%%%%%%%%%%%%%%%%%%%%%%%%%%%%%%%%%%%%%%%%%%%%%%%%%%%%%%%%%%%%%%%
\subsection{Stochastic Galerkin method}
The unknown coefficient functions in~(\ref{pce-appr}) can be 
determined approximately either by a stochastic collocation technique 
or a stochastic Galerkin method, see~\cite{xiu-hesthaven,xiu-book}.
In this paper, we investigate the stochastic Galerkin method. 
The approach yields a linear dynamical system 
\begin{equation} \label{galerkin}
\begin{array}{rcl}
\hat{E} \dot{\hat{v}}(t) & = & \hat{A} \hat{v}(t) + \hat{B} u(t) \\[1ex]
\hat{w}(t) & = & \hat{C} \hat{v}(t) \\
\end{array}
\end{equation}
for $t \ge 0$ of the larger dimension~$mn$ for the inner variables, 
where $m := | \mathcal{I} |$ denotes the cardinality of the index set. 
Therein, $\hat{v} = (\hat{v}_{i_1},\ldots,\hat{v}_{i_m}) \in \real^{mn}$ and 
$\hat{w} = (\hat{w}_{i_1},\ldots,\hat{w}_{i_m}) \in \real^m$ 
represent approximations to the exact 
coefficient functions~(\ref{coefficients}). 
The initial values are $\hat{v}(0)=0$ due to $x(0,p) = 0$ for all~$p$.
The constant matrices $\hat{A},\hat{E} \in \real^{mn \times mn}$, 
$\hat{B} \in \real^{mn}$ and $\hat{C} \in \real^{m \times mn}$ 
follow directly from the matrices in the dynamical system~(\ref{dae}) 
and the probability distribution of the parameters by integrals 
in the random space. 
For more details on the derivation of the system~(\ref{galerkin}), 
see~\cite[Sect.~2.4]{pulch-maten}, for example.

The linear dynamical system~(\ref{dae}) is assumed to be SISO, 
whereas the Galerkin system~(\ref{galerkin}) becomes 
single-input-multiple-output (SIMO).
The system~(\ref{galerkin}) may be unstable even though all 
original systems~(\ref{dae}) are stable, see~\cite{sonday}. 
However, this loss of stability appears rather seldom.
The larger dynamical system~(\ref{galerkin}) exhibits an own input-output 
behaviour with a transfer function 
$\hat{H} : (\complex \backslash \Theta) \rightarrow \complex^n$ 
for some finite set of poles $\Theta \subset \complex$ given by
\begin{equation} \label{transfer-galerkin}
\hat{H} (s) := 
\hat{C} \left( s \hat{E} - \hat{A} \right)^{-1} \hat{B} 
\qquad \mbox{for} \;\; s \in \complex \backslash \Theta . 
\end{equation}
In the vector $\hat{H} = (\hat{H}_1,\ldots,\hat{H}_m)^\top$, 
the components represent approximations to the coefficient functions of an 
orthogonal expansion for the original transfer function~(\ref{transferfcn})
in the frequency domain, 
see~\cite[Sect.~3.3]{pulch-maten-augustin}.

%%%%%%%%%%%%%%%%%%%%%%%%%%%%%%%%%%%%%%%%%%%%%%%%%%%%%%%%%%%%%%%%%%%%%%%%%%%%%%
\subsection{Model order reduction}
Since the system~(\ref{galerkin}) of the stochastic Galerkin method is 
huge for high numbers of random parameters and index sets 
like~(\ref{index-large}), it represents an excellent candidate for 
a model order reduction (MOR). 
A reduction of a stochastic Galerkin system was investigated 
in~\cite{mi,pulch-maten-augustin,pulch-maten,zou} 
using Krylov subspace methods and 
in~\cite{pulch-scee} using balanced truncation.
The task consists in the construction of a smaller dynamical system
\begin{equation} \label{galerkin-reduced}
\begin{array}{rcl}
\bar{E} \dot{\bar{v}}(t) & = & 
\bar{A} \bar{v}(t) + \bar{B} u(t) \\[1ex]
\bar{w}(t) & = & \bar{C} \bar{v}(t) \\
\end{array}
\end{equation}
with dimension~$r \ll nm$ of the inner variables and 
the same number of outputs~$m$. 
Hence the sizes of the matrices are 
$\bar{A},\bar{E} \in \real^{r \times r}$, 
$\bar{B} \in \real^r$, $\bar{C} \in \real^{m \times r}$. 
In projection-based model reduction, projection matrices 
$T_{\rm l}, T_{\rm r} \in \real^{mn \times r}$ are determined and the 
reduced matrices read as 
\begin{equation} \label{matrices-reduced}
\bar{A} = T_{\rm l}^\top \hat{A} T_{\rm r}, \quad
\bar{B} = T_{\rm l}^\top \hat{B}, \quad
\bar{C} = \hat{C} T_{\rm r}, \quad
\bar{E} = T_{\rm l}^\top \hat{E} T_{\rm r} . 
\end{equation}
The property $\hat{w}(t) \approx \bar{w}(t)$ for $t \ge 0$ 
is desired with respect to some norm in the time domain. 
The reduced system~(\ref{galerkin-reduced}) yields  
\begin{equation} \label{mor-output}
\bar{y}^{(\mathcal{I})} (t,p) = 
\sum_{i \in \mathcal{I}} \bar{w}_i(t) \Phi_i(p) 
\end{equation}
as an approximation to the random quantity of interest.

Let $\bar{H} = (\bar{H}_1,\ldots,\bar{H}_m)^\top$ be the transfer 
function of the reduced system~(\ref{galerkin-reduced}).
The difference between two transfer functions can be quantified by 
Hardy norms, see~\cite[Ch.~5.7]{antoulas}.
We apply these norms componentwise in the following.
%The $\mathcal{H}_2$-norm and $\mathcal{H}_\infty$-norm
%of the differences reads as
%$$ \begin{array}{lcl}
%\left\| \hat{H}_i - \bar{H}_i \right\|_{\mathcal{H}_2} & = & 
%\sqrt{ \displaystyle \frac{1}{2\pi} \int_{-\infty}^{+\infty} 
%\left| \hat{H}_i ({\rm i} \, \omega) - 
%\bar{H}_i ({\rm i} \, \omega) \right|^2 \; {\rm d}\omega } \\[3ex]
%\left\| \hat{H}_i - \bar{H}_i \right\|_{\mathcal{H}_{\infty}} & = & 
%\displaystyle \sup_{\omega \in \real} 
%\left| \hat{H}_i ({\rm i} \, \omega) - 
%\bar{H}_i ({\rm i} \, \omega) \right| \\
%\end{array} $$
%for $i=1,\ldots,M$.
In the time domain, the usual integral norm
$$ \| u \|_{\mathcal{L}^2[0,\infty)} = 
\sqrt{ \int_0^\infty u(t)^2 \; \mbox{d}t } $$
of the Hilbert space $\mathcal{L}^2[0,\infty)$ is used.
The next lemma recalls some definitions as well as general results 
on the input-output relations of linear dynamical systems.
The proof follows from the statements in~\cite[Ch.~2]{doyle-book}.
Associated error measures are also illustrated 
in~\cite[Sect.~2.4]{benner-gugercin-willcox}.

%%% Lemma: %%%%%%%%%%%%%%%%%%%%%%%%%%%%%%%%%%%%%%%%%%%%%%%%%%%%%%%%%%%%%%%%%
\begin{lemma} \label{lemma:estimates}
Let the initial values of a linear dynamical system be zero.
If $G$ is the transfer function of a stable and proper SISO system, 
then the $\mathcal{H}_\infty$-norm
\begin{equation} \label{hinf-norm-general}
\left\| G \right\|_{\mathcal{H}_{\infty}} = 
\sup_{\omega \in \real} \left| G ({\rm i} \, \omega) \right| 
\end{equation}
exists and the input-output exhibits the bound
\begin{equation} \label{l2-estimate} 
\| y \|_{\mathcal{L}^2 [0,\infty)} \le 
\left\| G \right\|_{\mathcal{H}_{\infty}} 
\| u \|_{\mathcal{L}^2 [0,\infty)} . 
\end{equation}
If $G$ is the transfer function of a stable and strictly proper SISO system, 
then the $\mathcal{H}_2$-norm
\begin{equation} \label{htwo-norm-general}
\left\| G \right\|_{\mathcal{H}_2} =
\sqrt{ \displaystyle \frac{1}{2\pi} \int_{-\infty}^{+\infty} 
\left| G ({\rm i} \, \omega) \right|^2 \; {\rm d}\omega } 
\end{equation}
exists and the input-output satisfies
\begin{equation} \label{max-estimate}
\sup_{t \ge 0} | y(t) | \le 
\left\| G \right\|_{\mathcal{H}_2} 
\| u \|_{\mathcal{L}^2 [0,\infty)} . 
\end{equation}
\end{lemma}
%%%%%%%%%%%%%%%%%%%%%%%%%%%%%%%%%%%%%%%%%%%%%%%%%%%%%%%%%%%%%%%%%%%%%%%%%%%%%
We will consider the $m$~outputs of the stochastic Galerkin 
system~(\ref{galerkin}) separately as SISO systems, 
which allows for using Lemma~\ref{lemma:estimates} in 
Section~\ref{sec:sparsebasis} and Section~\ref{sec:newbasis}.

%%%%%%%%%%%%%%%%%%%%%%%%%%%%%%%%%%%%%%%%%%%%%%%%%%%%%%%%%%%%%%%%%%%%%%%%%%%%%
%%%                    Sparsening of a basis                              %%%
%%%%%%%%%%%%%%%%%%%%%%%%%%%%%%%%%%%%%%%%%%%%%%%%%%%%%%%%%%%%%%%%%%%%%%%%%%%%%

\section{Sparsification of a given basis}
\label{sec:sparsebasis}
In this section, we suppose that a large canonical basis is given a priori, 
where the representation is shortened by neglecting basis functions. 

%%%%%%%%%%%%%%%%%%%%%%%%%%%%%%%%%%%%%%%%%%%%%%%%%%%%%%%%%%%%%%%%%%%%%%%%%%%%%
\subsection{Selection of a reduced basis}
\label{sec:sparsening1}
For the quantity of interest, 
the stochastic Galerkin method yields the approximation
\begin{equation} \label{galerkin-output}
\hat{y}^{(\mathcal{I})}(t,p) = \displaystyle
\sum_{i \in \mathcal{I}} \hat{w}_i(t) \Phi_i(p) . 
\end{equation}
%The resulting total error exhibits the upper bound
%$$ \begin{array}{rcl}
%\left\| y(t,\cdot)-\hat{y}^{(\mathcal{I})}(t,\cdot) \right\|_{\mathcal{L}^2} 
%& \le & 
%\left\| y(t,\cdot) - y^{(\mathcal{I})}(t,\cdot) \right\|_{\mathcal{L}^2} +
%\left\| y^{(\mathcal{I})} (t,\cdot) - 
%\hat{y}^{(\mathcal{I})}(t,\cdot) \right\|_{\mathcal{L}^2} \\[1ex] 
%& = & \displaystyle
%\sqrt{\sum_{i \notin \mathcal{I}} w_i(t)^2} + 
%\sqrt{\sum_{i \in \mathcal{I}} (w_i(t) - \hat{w}_i(t))^2} , \\
%\end{array} $$
%which consists of the truncation error and the error of the Galerkin approach.
%We assume that the set~$\mathcal{I}$ is chosen sufficiently large such 
%that both error parts are sufficiently small.
Now we choose a subset $\mathcal{I}' \subset \mathcal{I}$ or, 
equivalently, a selection of the outputs from the 
larger coupled system~(\ref{galerkin}).
The new approximation reads as 
\begin{equation} \label{sparsed-output}
\tilde{y}^{(\mathcal{I}')}(t,p) = \displaystyle
\sum_{i \in \mathcal{I}'} \hat{w}_i(t) \Phi_i(p) . 
\end{equation}
We denote this function by $\tilde{y}^{(\mathcal{I}')}$, 
because $\hat{y}^{(\mathcal{I}')}$ is devoted to the solution from the 
stochastic Galerkin method with respect to the index set~$\mathcal{I}'$.
The total error can be estimated by
$$ \begin{array}{rcl}
\left\| y(t,\cdot) - \tilde{y}^{(\mathcal{I}')}(t,\cdot) \right\|_{\mathcal{L}^2(\Omega)} 
& \le &  
\left\| y(t,\cdot) - y^{(\mathcal{I})}(t,\cdot) \right\|_{\mathcal{L}^2(\Omega)} \\[1ex]
& & 
+ \left\| y^{(\mathcal{I})} (t,\cdot) - 
\hat{y}^{(\mathcal{I})}(t,\cdot) \right\|_{\mathcal{L}^2(\Omega)} \\[1ex] 
& & + \left\| \hat{y}^{(\mathcal{I})} (t,\cdot) - 
\tilde{y}^{(\mathcal{I}')}(t,\cdot) \right\|_{\mathcal{L}^2(\Omega)} \\
\end{array} $$
for each $t \ge 0$, where the upper bound consists of three terms:  
the truncation error, the error of the Galerkin method and an 
additional sparsification error. 
We assume that the set $\mathcal{I}$ is chosen sufficiently large such that 
the truncation error as well as the Galerkin error are negligible.

The following theorem highlights the properties of the reduction 
in this section.

%%%%%%%%%%%%%%%%%%%%%%%%%%%%%%%%%%%%%%%%%%%%%%%%%%%%%%%%%%%%%%%%%%%%%%%%%%%%%
\begin{theorem} \label{thm1}
The error between the quantity of interest~(\ref{galerkin-output}) 
and its sparse approximation~(\ref{sparsed-output}) exhibits the estimates
\begin{eqnarray}  
\displaystyle \label{errorbounds1-1} 
\sup_{t \ge 0} \left\| \hat{y}^{(\mathcal{I})} (t,\cdot) - 
\tilde{y}^{(\mathcal{I}')}(t,\cdot) \right\|_{\mathcal{L}^2(\Omega)} & \le & 
\sqrt{ \displaystyle \sum_{i \in \mathcal{I} \backslash \mathcal{I}'} 
\left\| \hat{H}_i \right\|_{\mathcal{H}_2}^2 } \;\;
\| u \|_{\mathcal{L}^2[0,\infty)} \\
\displaystyle \label{errorbounds1-2}
\left\| \hat{y}^{(\mathcal{I})} - \tilde{y}^{(\mathcal{I}')}
\right\|_{\mathcal{L}^2(\Omega) \times \mathcal{L}^2[0,\infty)} & \le & 
\sqrt{ \displaystyle \sum_{i \in \mathcal{I} \backslash \mathcal{I}'} 
\left\| \hat{H}_i \right\|_{\mathcal{H}_\infty}^2 } \;\;
\| u \|_{\mathcal{L}^2[0,\infty)} 
\end{eqnarray} 
provided that the Galerkin system~(\ref{galerkin}) is stable as well as 
strictly proper or proper, respectively. 
\end{theorem}
%%%%%%%%%%%%%%%%%%%%%%%%%%%%%%%%%%%%%%%%%%%%%%%%%%%%%%%%%%%%%%%%%%%%%%%%%%%%%

Proof: \nopagebreak

The error of the sparse approximation satisfies the estimate
$$ \left\| \hat{y}^{(\mathcal{I})} (t,\cdot) - 
\tilde{y}^{(\mathcal{I}')}(t,\cdot) \right\|_{\mathcal{L}^2(\Omega)} = 
\sqrt{\sum_{i \in \mathcal{I} \backslash \mathcal{I}'} \hat{w}_i(t)^2} 
\quad \mbox{for each}\;\; t \ge 0 $$
due to the orthonormality of the basis functions. 

Firstly, using estimate~(\ref{max-estimate}) from Lemma~\ref{lemma:estimates}, 
we obtain 
$$ \sup_{t \ge 0} \left\| \hat{y}^{(\mathcal{I})} (t,\cdot) - 
\tilde{y}^{(\mathcal{I}')}(t,\cdot) \right\|_{\mathcal{L}^2(\Omega)}^2 \le 
\sum_{i \in \mathcal{I} \backslash \mathcal{I}'} 
\sup_{t \ge 0} \hat{w}_i(t)^2 \le
\sum_{i \in \mathcal{I} \backslash \mathcal{I}'} 
\left\| \hat{H}_i \right\|_{\mathcal{H}_2}^2
\| u \|_{\mathcal{L}^2[0,\infty)}^2 . $$
Taking the square root yields the statement.
 
Secondly, the estimate~(\ref{l2-estimate}) from Lemma~\ref{lemma:estimates}
implies
$$ \begin{array}{l} \displaystyle
\left\| \hat{y}^{(\mathcal{I})} - \tilde{y}^{(\mathcal{I}')}
\right\|_{\mathcal{L}^2(\Omega) \times \mathcal{L}^2[0,\infty)}^2 \; = \;
\displaystyle \int_0^\infty \int_{\Pi} 
\left( \hat{y}^{(\mathcal{I})} (t,p) - \tilde{y}^{(\mathcal{I}')} (t,p) 
\right)^2 \rho(p) \; \mbox{d}p \; \mbox{d}t \\[3ex]
= \; 
\displaystyle \int_0^\infty \sum_{i \in \mathcal{I} \backslash \mathcal{I}'} 
\hat{w}_i(t)^2 \; \mbox{d}t \; = \;
\sum_{i \in \mathcal{I} \backslash \mathcal{I}'} 
\| \hat{w}_i \|_{\mathcal{L}^2[0,\infty)}^2 \; \le \; 
\sum_{i \in \mathcal{I} \backslash \mathcal{I}'} 
\left\| \hat{H}_i \right\|_{\mathcal{H}_{\infty}}^2
\| u \|_{\mathcal{L}^2[0,\infty)}^2 . \\
\end{array} $$
Again the square root operation shows the claimed estimate. 
\hfill $\Box$

\medskip

The involved Hardy norms can be computed a priori from the matrices 
in the Galerkin system~(\ref{galerkin}).
On the one hand, the best choice of an index set $\mathcal{I}'$ with exactly 
$k$~elements for minimising the upper error bound 
(\ref{errorbounds1-1}) or (\ref{errorbounds1-2}) reads as
$$ \mathcal{I}' = {\rm arg} \, {\rm min} \left\{
\sum_{i \in \mathcal{I} \backslash \mathcal{I}'} 
\left\| \hat{H}_i \right\|^2 \; : \; 
| \mathcal{I}' | = k \right\} . $$
On the other hand, the smallest index set $\mathcal{I}'$, 
where the upper error bound (\ref{errorbounds1-1}) or (\ref{errorbounds1-2}) 
is below a predetermined threshold~$\delta > 0$ for a unit norm 
of the input becomes
$$ \mathcal{I}' = {\rm arg} \, {\rm min} \left\{ | \mathcal{I}' | \; : \;
\sum_{i \in \mathcal{I} \backslash \mathcal{I}'} 
\left\| \hat{H}_i \right\|^2 < \delta^2 \right\} . $$

A drawback of this approximation is that the reduction is based on an 
upper bound of the sparsification error, which is not sharp in general. 

A measure for the quality of the approximation in a relative sense 
(input has unit norm) is defined as follows. 
Let $\mathcal{I}_r' \subseteq \mathcal{I}$ be an index set associated 
to the $r$ largest norms.
Although this set is not always unique, the resulting sum of the largest 
norms is unique. 
We arrange the ratio
\begin{equation} \label{theta}
\theta_r := 
\left( \sum_{i \in \mathcal{I}_r'} \left\| \hat{H}_i \right\|^2 
\right)^\frac{1}{2} /
\left( \sum_{i \in \mathcal{I}} \left\| \hat{H}_i \right\|^2 
\right)^\frac{1}{2}
\quad \mbox{for}\;\; r = 1,\ldots,|\mathcal{I}| . 
\end{equation}
Obviously, it holds that 
$0 \le \theta_1 \le \theta_2 \le \theta_3 \le \cdots \le 
\theta_{|\mathcal{I}|} = 1$.

%%%%%%%%%%%%%%%%%%%%%%%%%%%%%%%%%%%%%%%%%%%%%%%%%%%%%%%%%%%%%%%%%%%%%%%%%%%%%
\subsection{Reduction of the stochastic Galerkin system}
\label{sec:galerkin-shortened}
The sparse approximation, which is described in Section~\ref{sec:sparsening1}, 
represents a simplification of the representation for 
the quantity of interest. 
If still the system~(\ref{galerkin}) from the stochastic Galerkin method 
is solved, then computational effort is not saved yet. 
An idea for a workload reduction is to repeat 
the stochastic Galerkin approach for 
the smaller index set~$\mathcal{I}'$ instead of 
the initial index set~$\mathcal{I}$.
Consequently, a smaller system~(\ref{galerkin}) appears with  
dimension $| \mathcal{I}'| N$ instead of $| \mathcal{I} | N$
and $| \mathcal{I}'|$ outputs instead of $| \mathcal{I} |$ outputs.
If the matrices of the system~(\ref{galerkin}) have already been computed 
for the set~$\mathcal{I}$, then we obtain the reduced system
just by discarding the rows and columns of the matrices, which do not 
belong to basis functions in the set~$\mathcal{I}'$.
However, this approach for a reduction is critical due to two reasons:
\begin{enumerate}
\item
In the Galerkin approach, 
the subspace spanned by $\{ \Phi_i \, : \, i \in \mathcal{I} \}$ 
is required for both the approximation to the exact solution and to keep 
the associated residual small. 
The best approximation to the solution may still be sufficiently accurate 
in the reduced basis, while the error of the Galerkin method can increase 
too much.
\item
The reduction relies on information about the output only. 
The inner variables could exhibit a different behaviour with respect to 
the stochastic modes. 
Thus crucial interactions of the inner variables may be removed by 
the downsizing of the system matrices.  
\end{enumerate}
Nevertheless, we apply this reduction straightforward and check 
the reduction error a posteriori by the Hardy norms of the difference 
in the transfer functions.
Removing rows and columns in the matrices of the 
system~(\ref{galerkin}) yields a system~(\ref{galerkin-reduced}) 
with matrices~(\ref{matrices-reduced}), where the projection matrix 
$T_{\rm l} = T_{\rm r}$ results from taking rows and columns out of 
a square identity matrix. 
Hence we obtain a special case of a projection-based MOR. 
An additional modification is required in the output matrix, 
since the number of outputs is assumed to be the identical 
in the systems~(\ref{galerkin}) and~(\ref{galerkin-reduced}).
Thus we do not delete rows but change those rows into zero vectors.
It follows that 
\begin{equation} \label{output-matrix-shortened}
\bar{C} = D \hat{C} T_{\rm r}
\end{equation} 
with a square diagonal matrix~$D$ consisting of zeros and ones. 
Now the quantity of interest~$\hat{y}^{(\mathcal{I}')}$ 
from the smaller Galerkin system~(\ref{galerkin}) is identical 
to $\bar{y}^{(\mathcal{I})}$ in~(\ref{mor-output}) 
from the reduced system~(\ref{galerkin-reduced}).

We specify bounds on the difference between the random 
quantity of interest resulting from general systems~(\ref{galerkin}) 
and~(\ref{galerkin-reduced}) in the following theorem. 
These estimates are also crucial within Section~\ref{sec:newbasis}.

%%%%%%%%%%%%%%%%%%%%%%%%%%%%%%%%%%%%%%%%%%%%%%%%%%%%%%%%%%%%%%%%%%%%%%%%%%%%%
\begin{theorem} \label{thm2}
The difference between the output from the stochastic Galerkin 
system~(\ref{galerkin}) and the output from the reduced 
system~(\ref{galerkin-reduced}) satisfies the estimates
\begin{eqnarray} 
\label{errorbounds2-h2} \displaystyle
\sup_{t \ge 0} \left\| \hat{y}^{(\mathcal{I})} (t,\cdot) - 
\bar{y}^{(\mathcal{I})}(t,\cdot) \right\|_{\mathcal{L}^2(\Omega)} & \le & 
\sqrt{ \displaystyle \sum_{i \in \mathcal{I}} 
\left\| \hat{H}_i - \bar{H}_i \right\|_{\mathcal{H}_2}^2 } \;\;
\| u \|_{\mathcal{L}^2[0,\infty)} \\
\label{errorbounds2-hinf} \displaystyle 
\left\| \hat{y}^{(\mathcal{I})} - \bar{y}^{(\mathcal{I})}
\right\|_{\mathcal{L}^2(\Omega) \times \mathcal{L}^2[0,\infty)} & \le & 
\sqrt{ \displaystyle \sum_{i \in \mathcal{I}} 
\left\| \hat{H}_i - \bar{H}_i \right\|_{\mathcal{H}_\infty}^2 } \;\;
\| u \|_{\mathcal{L}^2[0,\infty)} 
\end{eqnarray}
provided that the involved systems are stable and the norms are 
finite. 
\end{theorem}
%%%%%%%%%%%%%%%%%%%%%%%%%%%%%%%%%%%%%%%%%%%%%%%%%%%%%%%%%%%%%%%%%%%%%%%%%%%%%
The proof of these estimates is straightforward by repeating the 
concept of the proof for Theorem~\ref{thm1}. 

Since the smaller Galerkin system yields the quantity of interest
$\hat{y}^{(\mathcal{I}')} = \bar{y}^{(\mathcal{I})}$, 
it holds that $\bar{H}_i = 0$ for $i \notin \mathcal{I}'$ due to 
the definition of the output matrix~(\ref{output-matrix-shortened}).
Consequently, the upper error estimates~(\ref{errorbounds2-h2}) 
and~(\ref{errorbounds2-hinf}) are bounded from below by 
\begin{equation} \label{bound-below}
\sum_{i \in \mathcal{I}} 
\left\| \hat{H}_i - \bar{H}_i \right\| \ge 
\sum_{i \notin \mathcal{I}} 
\left\| \hat{H}_i \right\|^2 . 
\end{equation}
The approximation quality of the smaller Galerkin system is limited 
by the rate of decay of the Hardy norms in the original Galerkin system.

%%%%%%%%%%%%%%%%%%%%%%%%%%%%%%%%%%%%%%%%%%%%%%%%%%%%%%%%%%%%%%%%%%%%%%%%%%%%%
%%%                  Construction of a Small Basis                        %%%
%%%%%%%%%%%%%%%%%%%%%%%%%%%%%%%%%%%%%%%%%%%%%%%%%%%%%%%%%%%%%%%%%%%%%%%%%%%%%

\section{Construction of a small basis by MOR}
\label{sec:newbasis}
We discuss the derivation of an alternative basis for a sparse 
approximation now.

%%%%%%%%%%%%%%%%%%%%%%%%%%%%%%%%%%%%%%%%%%%%%%%%%%%%%%%%%%%%%%%%%%%%%%%%%%%%%
\subsection{Definition of basis functions}
\label{sec:newbasis-def}
The system~(\ref{galerkin-reduced}) from the MOR should feature a 
dimensionality $r \ll nm$, where $m=|\mathcal{I}|$ is the number of 
basis functions included in~(\ref{pce-appr}). 
In the case of a high-dimensional parameter space, 
the number~$m$ becomes huge. 
Often the MOR techniques are very efficient and thus produce a 
sufficiently accurate system with a dimensionality $r \ll m$.
In this case, we construct an alternative representation for the 
quantity of interest. 
Using the random process~(\ref{mor-output}) 
from the reduced system~(\ref{galerkin-reduced}), 
it holds that
$$ \bar{y}^{(\mathcal{I})}(t,p) = \displaystyle
\sum_{i \in \mathcal{I}} \bar{w}_i(t) \Phi_i(p) =
\sum_{i \in \mathcal{I}} 
\left[ \sum_{j=1}^r \bar{c}_{ij} \bar{v}_j(t) \right] \Phi_i(p) = 
\sum_{j=1}^r \bar{v}_j(t)
\left[ \sum_{i \in \mathcal{I}} \bar{c}_{ij}\Phi_i(p) \right] . $$
Now we define new functions on the parameter space~$\Pi$ by
\begin{equation} \label{newbasis}
\Psi_j(p) := \sum_{i \in \mathcal{I}} \bar{c}_{ij}\Phi_i(p) 
\quad \mbox{for}\;\; j=1,\ldots,r . 
\end{equation}
% These functions are not necessarily linear independent. 
It is rather unlikely that the functions~(\ref{newbasis}) 
are linearly dependent. 
However, the $r$ columns of the matrix~$\bar{C}$ often do not exhibit 
a full numerical rank with respect to a standard machine precision, 
even if the projection matrix~$T_{\rm r}$ in~(\ref{matrices-reduced}) 
has full (numerical) rank. 
This property allows for a small amount of an additional basis reduction 
to be analysed in the following subsection. 

From above, we obtain the sparse representation 
\begin{equation} \label{bary}
\bar{y}^{(\mathcal{I})}(t,p) = \displaystyle
\sum_{j=1}^r \bar{v}_j(t) \Psi_j(p) 
\end{equation}
with no sparsification error but only the error of the reduction.
The sparsity is measured by $\sigma = r/m$. 
The low-dimensional system~(\ref{galerkin-reduced}) yields both this 
sparse formulation and a cheap simulation in the time domain to 
obtain the coefficient functions $\bar{v}_1,\ldots,\bar{v}_r$. 
The bounds of Theorem~\ref{thm2} are valid for the difference between  
the sparse approximation~(\ref{bary}) from 
the system~(\ref{galerkin-reduced}) 
and the quantity of interest~(\ref{galerkin-output})
form the system~(\ref{galerkin}).

%%%%%%%%%%%%%%%%%%%%%%%%%%%%%%%%%%%%%%%%%%%%%%%%%%%%%%%%%%%%%%%%%%%%%%%%%%%%%
\subsection{Transformation to an orthonormal basis}
\label{sec:trafo-orthogonal}
The basis $\{ \Psi_1 , \ldots , \Psi_r \}$ from~(\ref{newbasis}) is not 
orthogonal. 
Nevertheless, we can transform this basis to an orthonormal description 
using the matrix~$\bar{C} \in \real^{m \times r}$ 
of the reduced system~(\ref{galerkin-reduced}) only.
Let $\Phi := (\Phi_1,\ldots,\Phi_m)^\top$ and 
$\Psi := (\Psi_1,\ldots,\Psi_r)^\top$ be column vectors.
It holds that $\Psi = \bar{C}^\top \Phi$ due to~(\ref{newbasis}). 

We apply a singular value decomposition (SVD), 
see~\cite[Ch.~2.5]{golub-loan}, as a tool for two purposes: 
to construct an orthonormal basis and to remove unessential parts 
associated to a numerical rank deficiency of the output matrix.
The SVD of the output matrix reads as 
\begin{equation} \label{svd}
\bar{C} = U S Q \quad \mbox{with} \quad 
S = \begin{pmatrix} S' \\ 0 \\ \end{pmatrix} 
\end{equation} 
including orthogonal matrices $U \in \real^{m \times m}$, 
$Q \in \real^{r \times r}$ and a matrix $S \in \real^{m \times r}$, 
whose part $S' \in \real^{r \times r}$ is a diagonal matrix with 
the singular values $s_1 \ge s_2 \ge \cdots \ge s_r > 0$. 

It is straightforward to show the equality
\begin{equation} \label{psistar} 
\Psi^* := {S'}^{-1} Q \Psi = (I,0) U^\top \Phi 
\end{equation}
with the identity matrix $I \in \real^{r \times r}$.
Hence we have identified an orthonormal basis 
$\Psi^* = (\Psi_1^* , \ldots , \Psi_r^*)^\top$ satisfying
${\rm span} \{ \Psi_1^* , \ldots , \Psi_r^* \} = 
{\rm span} \{ \Psi_1 , \ldots , \Psi_r \}$. 
On the one hand, the basis functions~$\Psi^*$ can be calculated from 
the original basis $\Phi$ via~(\ref{psistar}) using the matrix~$U$.
On the other hand, it holds that
\begin{equation} \label{output2}
\bar{y}^{(\mathcal{I})}(t,p) = \displaystyle
\sum_{j=1}^r \bar{v}_j(t) \Psi_j(p) = 
\sum_{\ell=1}^r \left[ \sum_{j = 1}^r s_{\ell} q_{\ell j} \bar{v}_j(t) \right] 
\Psi_{\ell}^* (p) =:
\sum_{\ell=1}^r \bar{v}_{\ell}^*(t) \Psi_{\ell}^*(p) 
\end{equation}
with new coefficient functions $\bar{v}_{\ell}^*$
computable from original coefficient functions~$\bar{v}_j$ 
by the matrix $Q=(q_{ij})$ due to $\Psi = Q^\top S' \Psi^*$. 
Hence an explicit calculation of the basis functions~$\Psi$ is 
never required.  
Now we define an approximation by a truncation of the 
finite sum in~(\ref{output2}), i.e., 
\begin{equation} \label{output-deflation}
\breve{y}^{(r')}(t,p) := \displaystyle \sum_{\ell=1}^{r'} 
\left[ \sum_{j = 1}^r s_{\ell} q_{\ell j} \bar{v}_j(t) \right] 
\Psi_{\ell}^* (p)
\end{equation}
for an $r' \in \{ 1,\ldots,r-1 \}$. 
In our context, the following quality of the approximation is guaranteed.
%%% Lemma: %%%%%%%%%%%%%%%%%%%%%%%%%%%%%%%%%%%%%%%%%%%%%%%%%%%%%%%%%%%%%%%%%
\begin{theorem} \label{thm:deflation}
The output~(\ref{output2}) and its approximation~(\ref{output-deflation}) 
satisfy the error estimates
\begin{eqnarray} 
\label{deflation1} 
\left\| \bar{y}^{(\mathcal{I})} (t,\cdot) - 
\breve{y}^{(r')} (t,\cdot) \right\|_{\mathcal{L}^2(\Omega)} & \le &  
\sqrt{r - r'} \; s_{r'+1} \; \sqrt{ \sum_{j=1}^r \bar{v}_j(t)^2 }
\quad \mbox{for}\;\; t\ge 0 \quad \mbox{} \\
\label{deflation2} 
\left\| \bar{y}^{(\mathcal{I})} - \breve{y}^{(r')}
\right\|_{\mathcal{L}^2(\Omega) \times \mathcal{L}^2[0,\infty)} & \le &  
\sqrt{r - r'} \; s_{r'+1} \; 
\sqrt{ \sum_{j=1}^r \| \bar{v}_j \|_{\mathcal{L}^2[0,\infty)}^2 }
\end{eqnarray}
with $r'=1,\ldots,r-1$.
\end{theorem}
%%%%%%%%%%%%%%%%%%%%%%%%%%%%%%%%%%%%%%%%%%%%%%%%%%%%%%%%%%%%%%%%%%%%%%%%%%%%%

Proof:

We obtain 
$$ \bar{y}^{(\mathcal{I})} (t,p) - 
\breve{y}^{(r')} (t,p) = 
\sum_{\ell = r'+1}^r 
\left[ \sum_{j=1}^r s_{\ell} q_{\ell j} \bar{v}_j(t) \right] 
\Psi_{\ell}^*(p) $$
for each $t \ge 0$ and $p \in \Pi$.
The orthonormality of the basis functions implies 
$$ \left\| \bar{y}^{(\mathcal{I})} (t,\cdot) - \breve{y}^{(r')} (t,\cdot) 
\right\|_{\mathcal{L}^2(\Omega)}^2 = 
\sum_{\ell = r'+1}^r 
\left[ \sum_{j=1}^r s_{\ell} q_{\ell j} \bar{v}_j(t) \right]^2 . $$
The Cauchy-Schwarz inequality and the orthogonality of the matrix~$Q$ yields
$$ \left( \sum_{j=1}^r s_{\ell} q_{\ell j} \bar{v}_j(t) \right)^2 \le 
s_{\ell}^2 \left( \sum_{j=1}^r q_{\ell j}^2 \right) 
\left( \sum_{j=1}^r \bar{v}_j(t)^2 \right) \le 
s_{r'+1}^2 \sum_{j=1}^r \bar{v}_j(t)^2 $$
uniformly for $\ell = r'+1 , \ldots , r$. 
Taking the sum of all $\ell = r'+1 , \ldots , r$ shows the 
estimate~(\ref{deflation1}). 
Employing the integral on the time domain~$[0,\infty)$ confirms 
the estimate~(\ref{deflation2}). 
\hfill $\Box$

\medskip

The time-dependent part of the upper bounds 
in~(\ref{deflation1}),(\ref{deflation2}) 
cannot be estimated further without additional assumptions. 
However, we apply this reduction step only to remove a 
numerical rank deficiency, where $r'$ is selected such that 
the dominant singular value $s_{r'+1}$ is tiny but still 
significantly above the machine precision.

Only the first~$r$ columns of the matrix~$U=(u_{ij})$ are 
applied to determine the new orthonormal basis in~(\ref{psistar}). 
Hence these columns provide a measure for the influence of an 
original basis function~$\Phi_i$ within the new basis. 
We define the values
\begin{equation} \label{kappa}
\kappa_i := \sqrt{ \sum_{j=1}^r u_{ij}^2 } 
\qquad \mbox{for}\;\; i=1,\ldots,m . 
\end{equation}
Since the matrix~$U$ is orthogonal, it holds that 
$0 \le \kappa_i \le 1$ for all~$i$ and $r \le m$.  
This measure of is independent of the quantification by Hardy norms 
as done in Section~\ref{sec:sparsebasis}.

%%%%%%%%%%%%%%%%%%%%%%%%%%%%%%%%%%%%%%%%%%%%%%%%%%%%%%%%%%%%%%%%%%%%%%%%%%%%%
%%%                         Test Example                                  %%%
%%%%%%%%%%%%%%%%%%%%%%%%%%%%%%%%%%%%%%%%%%%%%%%%%%%%%%%%%%%%%%%%%%%%%%%%%%%%%

\section{Illustrative example}
\label{sec:example}
Now we apply the reduction approaches from the previous sections 
to a test example. 
All numerical calculations are performed within the software package 
MATLAB (version R2014b), where the machine precision is around 
$\varepsilon_0 = 2 \cdot 10^{-16}$.

%%%%%%%%%%%%%%%%%%%%%%%%%%%%%%%%%%%%%%%%%%%%%%%%%%%%%%%%%%%%%%%%%%%%%%%%%%%%%
\subsection{Modelling of a low pass filter}
\begin{figure}[t]
\begin{center}
\includegraphics[width=14cm]{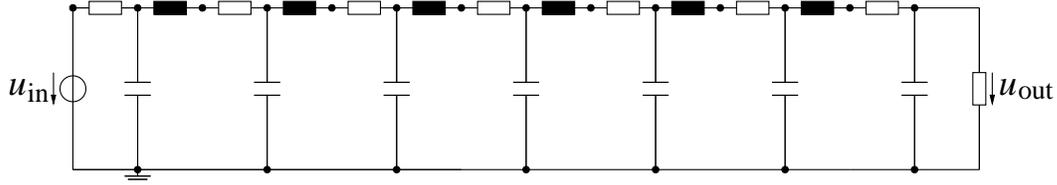} 
\end{center}
\caption{Electric circuit of a low pass filter.}
\label{fig:lowpassfilter}
\end{figure}

\begin{figure}
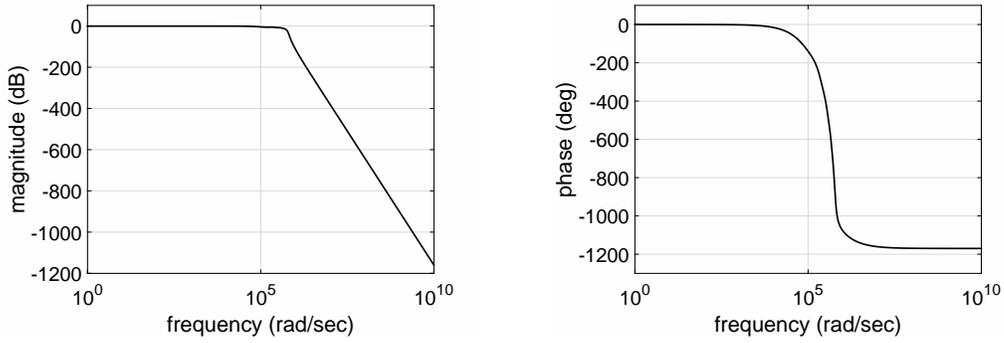

\begin{center}
\includegraphics[width=6cm]{bode_lowpass_1.eps}
\hspace{10mm}
\includegraphics[width=6cm]{bode_lowpass_2.eps}
\end{center}
\caption{Bode plot of the low pass filter for constant physical parameters: 
magnitude (left) and phase (right) of the transfer function.}
\label{fig:bode}
\end{figure}

As test example, we investigate the electric circuit of a low pass filter, 
see Figure~\ref{fig:lowpassfilter}. 
The electric circuit includes 7~capacitances, 6~inductances and 
8~conductances. 
Thus $q=21$ physical parameters occur. 
The system is SISO, because a single input voltage is supplied and the 
output voltage drops at a load conductor. 
A mathematical modelling yields a system~(\ref{dae}) of DAEs 
for the 14~node voltages and the 6 currents through the inductances. 
Hence the dimension of the inner variables is $n=20$. 
The (nilpotency) index of the DAE system results to one. 
Furthermore, the linear dynamical system is stable as well as 
strictly proper.
Figure~\ref{fig:bode} depicts the bode plot of this linear time-invariant 
system in the case of a constant choice for the parameters. 
The magnitude of the transfer function demonstrates the characteristics 
of a low pass filter.

%%% Figure: RLC circuit %%%%%%%%%%%%%%%%%%%%%%%%%%%%%%%%%%%%%%%%%%%%%%%%%%%%%
%%%%%%%%%%%%%%%%%%%%%%%%%%%%%%%%%%%%%%%%%%%%%%%%%%%%%%%%%%%%%%%%%%%%%%%%%%%%%

In the stochastic modelling, we replace all physical parameters 
by uniformly distributed random variables with a range of 10\% 
around their mean values. 
The mean values are the constant choice of the parameters from above.  
In the orthogonal expansion, we include all multivariate polynomials 
up to total degree $d=3$, which results in $m=2024$ basis functions. 
The stochastic Galerkin method generates a linear dynamical 
system~(\ref{galerkin}) of dimension $mn = 40480$.
We compute the involved matrices numerically, where the probabilistic 
integrals are approximated by a quadrature on a sparse grid with 
13329 nodes in the domain~$\Pi$. 
The output matrix~$\hat{C}$ is obtained directly, since the matrix~$C$ 
in~(\ref{dae}) does not depend on the physical parameters. 
Moreover, $\hat{C}$ exhibits full (numerical) rank.
Numerical computations confirm that the resulting linear dynamical 
system~(\ref{galerkin}) is stable as well as strictly proper. 

%%% Figure: Bode plot %%%%%%%%%%%%%%%%%%%%%%%%%%%%%%%%%%%%%%%%%%%%%%%%%%%%%%%
%%%%%%%%%%%%%%%%%%%%%%%%%%%%%%%%%%%%%%%%%%%%%%%%%%%%%%%%%%%%%%%%%%%%%%%%%%%%%

%%%%%%%%%%%%%%%%%%%%%%%%%%%%%%%%%%%%%%%%%%%%%%%%%%%%%%%%%%%%%%%%%%%%%%%%%%%%%
\subsection{Sparsification of the given basis}
\begin{figure}
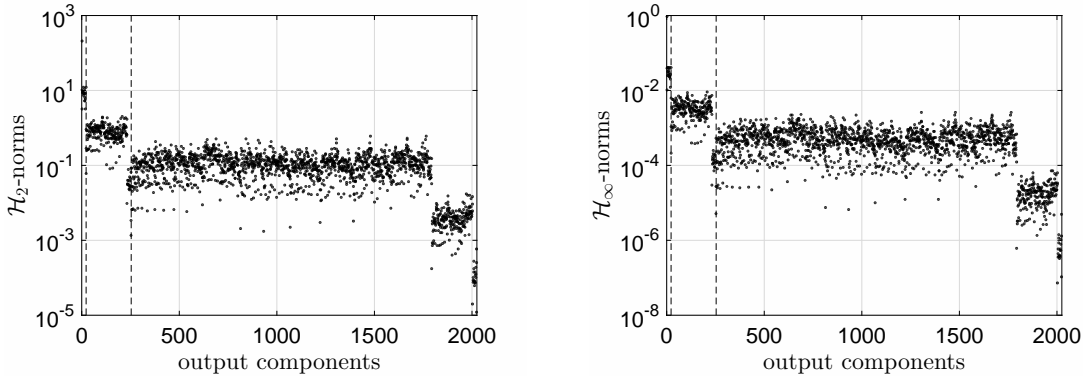

\begin{center}
\includegraphics[width=6.5cm]{htwo_gal.eps}
\hspace{10mm}
\includegraphics[width=6.5cm]{hinf_gal.eps}
\end{center}
\caption{$\mathcal{H}_2$-norms (left) and $\mathcal{H}_{\infty}$-norms 
(right) for the components of the transfer function in the stochastic 
Galerkin system. 
(The two dashed lines separate the components for polynomials of 
degree zero/one, degree two and degree three.)}
\label{fig:hnormsgal}
\end{figure}

\begin{figure}
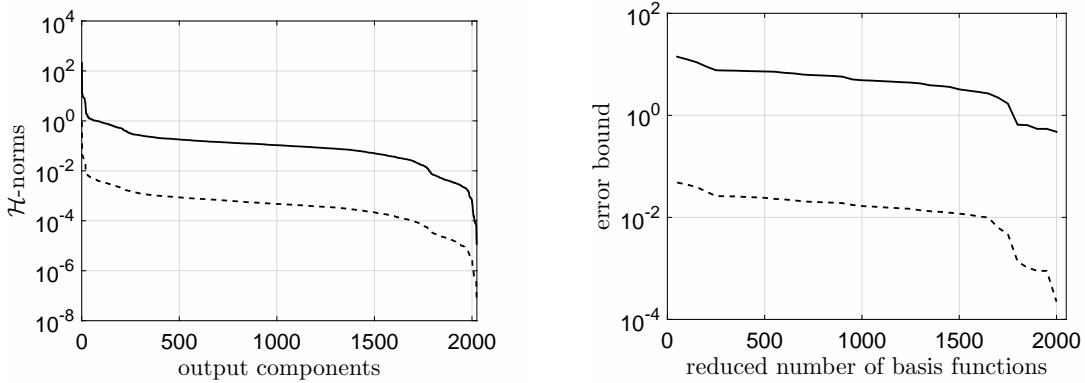

\begin{center}
\includegraphics[width=6.5cm]{hnorms_gal.eps}
\hspace{10mm}
\includegraphics[width=6.5cm,height=5.1cm]{error_reduce_gal.eps}
\end{center}
\caption{Hardy norms of the stochastic Galerkin system in descending order 
(left) and error bounds~(\ref{errorbounds2-h2}),(\ref{errorbounds2-hinf}) 
for reduced Galerkin system (right). 
Solid lines and dashed lines show data for $\mathcal{H}_2$-norms 
and $\mathcal{H}_{\infty}$-norms, respectively.}
\label{fig:hnormsgal2}
\end{figure}

We compute the $\mathcal{H}_2$-norms as well as $\mathcal{H}_{\infty}$-norms 
for the separate components of the transfer function associated to 
the stochastic Galerkin system~(\ref{galerkin}). 
For this purpose, the transfer function is evaluated on
a logarithmically spaced grid $0 \le \omega_1 < \omega_2 < \cdots < \omega_k$ 
inside the imaginary axis $\{ {\rm i}\omega \, : \, \omega \in \real \}$ 
due to symmetry reasons. 
For the $\mathcal{H}_{\infty}$-norm~(\ref{hinf-norm-general}), 
just the maximum of the finite set of absolute values is chosen 
as approximation. 
For the $\mathcal{H}_2$-norm~(\ref{htwo-norm-general}), 
the trapezoidal rule yields an approximation to the integral.
In each frequency point~$s={\rm i}\omega$, 
the computational effort for an evaluation of~(\ref{transfer-galerkin}) 
is dominated by the solution of a linear system of algebraic equations 
with coefficient matrix $s\hat{E}-\hat{A}$ and right-hand side~$\hat{B}$. 
Although the number of outputs is large, the computational work of 
the matrix-vector-multiplication with the sparse matrix~$\hat{C}$ is 
negligible.
The resulting norms are depicted for all output components in 
Figure~\ref{fig:hnormsgal}. 
Furthermore, Figure~\ref{fig:hnormsgal2} (left) shows these 
norms in a descending sequence.
We observe different orders of magnitudes for the norms, 
which mainly depend on the degree of the associated basis polynomials. 
Nevertheless, there are several components for degree two as well as 
degree three, whose norms have the same magnitude. 

\begin{figure}
\begin{center}
\includegraphics[width=6.5cm]{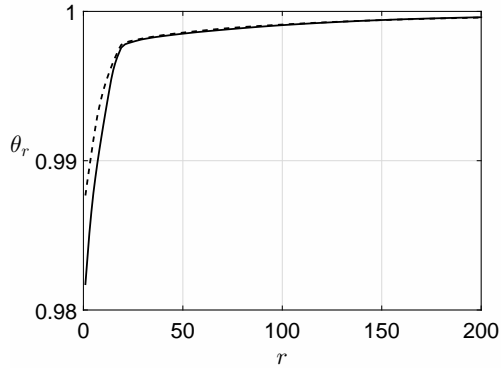}
\end{center}
\caption{Ratios~$\theta_r$ from~(\ref{theta}) for different 
index sets~$\mathcal{I}_r'$. 
Solid lines and dashed lines show data for $\mathcal{H}_2$-norms 
and $\mathcal{H}_{\infty}$-norms, respectively.}
\label{fig:theta}
\end{figure}

These results indicate some potential for a sparse
approximation~(\ref{sparsed-output}) as explained 
in Section~\ref{sec:sparsening1}.
On the one hand, the upper bounds from Theorem~\ref{thm1} could be 
evaluated for some selections of index sets. 
We omit the presentation for shortness.
On the other hand, Figure~\ref{fig:theta} illustrates the ratios $\theta_r$ 
from~(\ref{theta}) for index sets~$\mathcal{I}_r'$ with $r=1,\ldots,200$. 
We recognise that the values $\theta_r$ tend to one. 
Furthermore, Table~\ref{tab:dimensions1} shows the minimum 
cardinality~$r = |\mathcal{I}_r'|$ and the accompanying ratio $r/m$, 
which are necessary to achieve a value $\theta_r \ge 1 - \delta$
for several thresholds $\delta > 0$.

%%% Figure: H-norms %%%%%%%%%%%%%%%%%%%%%%%%%%%%%%%%%%%%%%%%%%%%%%%%%%%%%%%%%
%%%%%%%%%%%%%%%%%%%%%%%%%%%%%%%%%%%%%%%%%%%%%%%%%%%%%%%%%%%%%%%%%%%%%%%%%%%%%

%%% Figure: H-norms %%%%%%%%%%%%%%%%%%%%%%%%%%%%%%%%%%%%%%%%%%%%%%%%%%%%%%%%%
%%%%%%%%%%%%%%%%%%%%%%%%%%%%%%%%%%%%%%%%%%%%%%%%%%%%%%%%%%%%%%%%%%%%%%%%%%%%%

%%% Figure: Values theta %%%%%%%%%%%%%%%%%%%%%%%%%%%%%%%%%%%%%%%%%%%%%%%%%%%%
%%%%%%%%%%%%%%%%%%%%%%%%%%%%%%%%%%%%%%%%%%%%%%%%%%%%%%%%%%%%%%%%%%%%%%%%%%%%%

%%% Table: Minimum dimensions %%%%%%%%%%%%%%%%%%%%%%%%%%%%%%%%%%%%%%%%%%%%%%%
\begin{table}
\caption{Minimum cardinality~$r$ required for a 
value~$\theta_r \ge 1-\delta$ in~(\ref{theta}) and 
their ratios $r/m$ in percentage.}
\label{tab:dimensions1}
\begin{center}
\begin{tabular}{rccccc}
& & $\delta = 10^{-2}$ & $\delta = 10^{-3}$ & $\delta = 10^{-4}$ & 
$\delta = 10^{-5}$ \\ \hline
$\mathcal{H}_2$-norms & $r$ & 8 & 92 & 672 & 1350 \\ 
& $r/m$ & 0.4\% & 4.6\% & 33.2\% & 66.7\% \\ \hline
$\mathcal{H}_{\infty}$-norms & $r$ & 4 & 86 & 678 & 1331 \\
& $r/m$ & 0.2\% & 4.3\% & 33.5\% & 65.8\% \\ 
\end{tabular}
\end{center}
\end{table}
%%%%%%%%%%%%%%%%%%%%%%%%%%%%%%%%%%%%%%%%%%%%%%%%%%%%%%%%%%%%%%%%%%%%%%%%%%%%%

We also use the reduction technique from 
Section~\ref{sec:galerkin-shortened}, 
where the stochastic Galerkin system is shortened. 
In Table~\ref{tab:dimensions1}, the ratios also indicate the 
reduction of the dimensionality for the Galerkin system.
Furthermore, we perform this reduction for a sequence of 
index sets~$\mathcal{I}_r'$ with $10 \le r \le 2000$. 
The upper error bounds~(\ref{errorbounds2-h2}),(\ref{errorbounds2-hinf})
from Theorem~\ref{thm2} 
are evaluated assuming an input with unit norm, 
where the Hardy norms are approximated by the same approach as above. 
Figure~\ref{fig:hnormsgal2} (right) depicts the results. 
In agreement to Figure~\ref{fig:hnormsgal2} (left), we observe that 
the error estimates are bounded from below by~(\ref{bound-below}).
Hence only a moderate potential for a sparse approximation and 
a reduction of the stochastic Galerkin system occurs.

%%%%%%%%%%%%%%%%%%%%%%%%%%%%%%%%%%%%%%%%%%%%%%%%%%%%%%%%%%%%%%%%%%%%%%%%%%%%%
\subsection{New basis from MOR}
Now the approach from Section~\ref{sec:newbasis-def} is employed. 
We use a moment matching technique with a single expansion point 
in the frequency domain, see~\cite[Sect.~3.4]{freund}.
A Krylov subspace is defined by the input vector~$\hat{B}$. 
In our MOR technique, the output matrix~$\hat{C}$ is not involved 
and thus the computational effort is independent of the number of outputs. 
The Arnoldi algorithm yields a projection matrix $T_{\rm r}$, 
whose columns form an orthonormal basis.
We simply choose $T_{\rm l} = T_{\rm r}$. 
In the system~(\ref{galerkin-reduced}),
the smaller matrices follow from~(\ref{matrices-reduced}). 
We tried several choices for the expansion point on the real axis. 
The best instance resulted to $s = 5 \cdot 10^5$.
We can select an arbitrary dimension~$r \le nm$ of the 
reduced system~(\ref{galerkin-reduced}) by taking the first~$r$ 
vectors of the Arnoldi process.
Numerical computations confirm that the reduced 
systems~(\ref{galerkin-reduced}) are stable for $r=52,\ldots,100$
and instable for some $r \le 51$.
% For r=100,...200, all stable except r=137,138.
Moreover, numerical examinations show that the systems~(\ref{galerkin-reduced})
are strictly proper in all observed cases.
% Plots of transfer functions observed for r=10,...,100.

\begin{figure}
\begin{center}
\includegraphics[width=6.5cm]{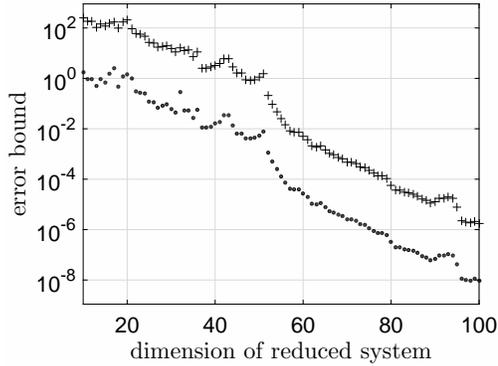} 
\end{center}
\caption{Error bounds from~(\ref{errorbounds2-h2})
with $\mathcal{H}_2$-norms (+) and~(\ref{errorbounds2-hinf}) 
with $\mathcal{H}_{\infty}$-norms ($\cdot$) 
for different dimensions of the reduced system.}
\label{fig:norm-difference}
\end{figure}

Figure~\ref{fig:norm-difference} illustrates the upper error bounds 
(\ref{errorbounds2-h2}),(\ref{errorbounds2-hinf}) for a unit norm 
of the input in the case of dimensions $r=10,\ldots,100$. 
The minimum dimensions, which are required to achieve 
error estimates below some thresholds, are given in 
Table~\ref{tab:dimensions2}.
We observe a rapid decay of the error bounds, which 
confirms a high potential for a reduction of the dimensionality 
in the linear dynamical system~(\ref{galerkin-reduced}).
Other MOR techniques, which involve the output matrix~$\hat{C}$, 
most likely generate even lower-dimensional systems~(\ref{galerkin-reduced}) 
of the same accuracy.

%%% Figure: Norms of differences %%%%%%%%%%%%%%%%%%%%%%%%%%%%%%%%%%%%%%%%%%%%
%%%%%%%%%%%%%%%%%%%%%%%%%%%%%%%%%%%%%%%%%%%%%%%%%%%%%%%%%%%%%%%%%%%%%%%%%%%%%

%%% Table: Minimum dimensions %%%%%%%%%%%%%%%%%%%%%%%%%%%%%%%%%%%%%%%%%%%%%%%
\begin{table}
\caption{Minimum dimensions required for an error bound below a 
threshold~$\delta$.}
\label{tab:dimensions2}
\begin{center}
\begin{tabular}{ccccc}
& $\delta = 10^{-1}$ & $\delta = 10^{-2}$ & $\delta = 10^{-3}$ & 
$\delta = 10^{-4}$ \\ \hline
bound~(\ref{errorbounds2-h2}) & 53 & 57 & 67 & 80 \\
bound~(\ref{errorbounds2-hinf}) & 27 & 45 & 53 & 56 
\end{tabular}
\end{center}
\end{table}
%%%%%%%%%%%%%%%%%%%%%%%%%%%%%%%%%%%%%%%%%%%%%%%%%%%%%%%%%%%%%%%%%%%%%%%%%%%%%

\begin{figure}
\begin{center}
\includegraphics[width=6.5cm]{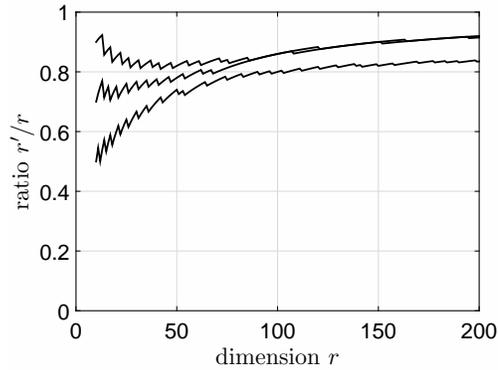} 
\end{center}
\caption{Ratios between the dimension~$r$ of the reduced systems 
and the dimension~$r'$ after a further reduction by the omission of 
all singular values below the thresholds $10^{-4}$ (lower line), 
$10^{-8}$ (center line) and $10^{-12}$ (upper line).}
\label{fig:deflation}
\end{figure}

We also use the technique described in Section~\ref{sec:trafo-orthogonal}, 
where the SVD of the output matrix is computed. 
An additional reduction from the number~$r$ in~(\ref{output2}) 
to a lower number~$r'$ motivated by Theorem~\ref{thm:deflation}
is feasible by neglecting all singular values below some threshold. 
Figure~\ref{fig:deflation} shows the ratios $r' / r$ for different 
dimensions~$r$ and three choices of the threshold. 
We recognise some gain in efficiency by removing the unessential parts. 
For example, a dimension $r=100$ can be decreased further to 
$r' = \mbox{80-90}$ depending on the required accuracy. 

%%% Figure: Deflation %%%%%%%%%%%%%%%%%%%%%%%%%%%%%%%%%%%%%%%%%%%%%%%%%%%%%%%
%%%%%%%%%%%%%%%%%%%%%%%%%%%%%%%%%%%%%%%%%%%%%%%%%%%%%%%%%%%%%%%%%%%%%%%%%%%%%

%%% Figure: Values Kappa %%%%%%%%%%%%%%%%%%%%%%%%%%%%%%%%%%%%%%%%%%%%%%%%%%%%
\begin{figure}
\begin{center}
\includegraphics[width=6.5cm]{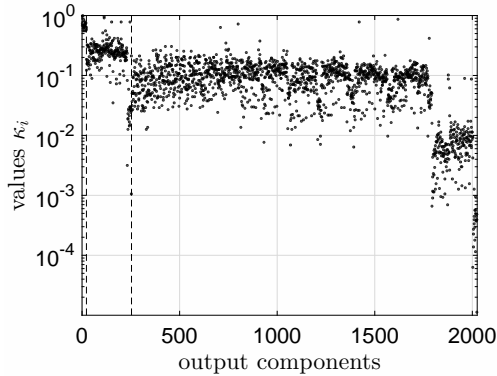} 
\end{center}
\caption{Values~$\kappa_i$ from~(\ref{kappa}) for all output components in 
the case of a reduced system with dimension $r=50$. 
(The two dashed lines separate the components for polynomials of 
degree zero/one, degree two and degree three.)}
\label{fig:kappa}
\end{figure}
%%%%%%%%%%%%%%%%%%%%%%%%%%%%%%%%%%%%%%%%%%%%%%%%%%%%%%%%%%%%%%%%%%%%%%%%%%%%%

Finally, we compute the values~$\kappa_i$ from~(\ref{kappa}) 
using the SVD~(\ref{svd}), 
where the reduced system exhibits the dimension $r=50$. 
The results are depicted in Figure~\ref{fig:kappa}. 
We recover the structure of the Hardy norms for the original stochastic  
Galerkin system~(\ref{galerkin}), cf. Figure~\ref{fig:hnormsgal}, 
although the computation of the values~(\ref{kappa}) is done 
completely different. 
Yet some components for basis functions of degree two and three 
feature relatively high numbers (close to the upper bound one) now. 
This property indicates that some components, which are less important
in the system~(\ref{galerkin}), become crucial for achieving our 
sparse approximation. 

%%%%%%%%%%%%%%%%%%%%%%%%%%%%%%%%%%%%%%%%%%%%%%%%%%%%%%%%%%%%%%%%%%%%%%%%%%%%%
%%%                          Conclusions                                  %%%
%%%%%%%%%%%%%%%%%%%%%%%%%%%%%%%%%%%%%%%%%%%%%%%%%%%%%%%%%%%%%%%%%%%%%%%%%%%%%

\section{Conclusions and outlook}
We examined two numerical techniques for the identification of a 
sparse representation, which approximates a quantity of interest from 
a random linear dynamical system. 
On the one hand, terms of an orthogonal expansion were omitted 
if their accompanying Hardy norms are relatively small. 
On the other hand, a projection-based model order reduction yields 
an alternative orthogonal expansion with a low number of basis functions. 
We proved upper error bounds for the sparse approximations in 
both techniques. 
In addition, we carried out numerical computations for a test example with 
high dimensionality. 
The results indicate that the second approach is more efficient than the 
first method. 
If the same accuracy is required with respect to the error bounds, 
then the second technique yields a sparse approximation with a 
lower number of basis functions. 
%Furthermore, we showed that the model order reduction yields an 
%alternative measure for the importance of the original basis functions, 
%where a singular value decomposition is applied to the output matrix. 
%In comparison to the quantification by Hardy norms, this measure 
%is cheaper to compute. 

A topic for further research is the construction of 
a sparse approximation in the case of random nonlinear dynamical systems. 
Model order reduction becomes harder in the nonlinear situation. 
For example, a strict estimation of approximation errors is critical 
or even impossible. 
As a simplification, one could investigate a linear dynamical system, 
whose quantity of interest depends nonlinearly on the state variables 
or the inner variables, respectively.

%%%%%%%%%%%%%%%%%%%%%%%%%%%%%%%%%%%%%%%%%%%%%%%%%%%%%%%%%%%%%%%%%%%%%%%%%%%%%
%%%                           References                                  %%%
%%%%%%%%%%%%%%%%%%%%%%%%%%%%%%%%%%%%%%%%%%%%%%%%%%%%%%%%%%%%%%%%%%%%%%%%%%%%%

\end{document}